\documentstyle[mydefs,12pt,psfig]{article}
\def\frak{\bs}

\newcommand{\be}{\begin{equation}}
\newcommand{\ee}{\end{equation}}
\newcommand{\la}{\left|\begin{array}{ccc}}
\newcommand{\ra}{\end{array}\right|}
\newcommand{\ma}{\frak}

\newcommand{\unity}{{ \rm \setlength{\unitlength}{1em}
\advance\textheight by 3cm
\advance\topmargin by -1.8cm
\advance\textwidth by 0.5cm
\begin{picture}(0.75,1)
\put(0,0){1}\put(0.34,0){\line(0,1){0.65}}
\end{picture} }}
\newcommand{\eins}{{\rm \setlength{\unitlength}{1em}
\begin{picture}(0.75,1)
\put(0,0){1}\put(0.34,0){\line(0,1){0.65}}
\end{picture}}}
\begin{document}

\input amssym.def
\input amssym.tex

\begin{center}
{ \LARGE \bf Spin Topological Quantum\\[0.2cm]
 Field Theories }

 \vspace*{1.5cm}

\begin{tabular}{c}
Anna Beliakova \\
%
 Institut de Recherche Math\'ematique Avanc\'ee
\\
Universit\'e Louis Pasteur, Strasbourg 
\end{tabular}
\vspace*{0.5cm}
\\August 96 \vspace*{0.5cm}

\end{center}

\parbox{10cm}{{\footnotesize {\bf Abstract}: 
Starting from the quantum group $U_q({\ma sl}(2,\Bbb C))$, we construct
operator invariants of
3-cobordisms with spin
structure, satisfying the requirements of
a topological quantum field theory and 
refining the Reshetikhin--Turaev and Turaev--Viro models.
We establish the relationship between these two refined models.}}

\section{Introduction}\label{ein}
This paper  is devoted to the refinement of the quantum invariants of 
3-manifolds taking into account spin structures.
The invariants of
Reshetikhin--Turaev type,  corresponding to the
quantum group  $U_q ({\ma sl}(2,\Bbb C))$ 
and  determined by
 a spin structure 
on a
closed 3-manifold, were first constructed by Blanchet [Bl],
 Kirby--Melvin [KM]  and Turaev [Tu].  
The idea of the construction was the following: 
Using a presentation of a closed 3-manifold $M$
by  surgery along a link $L$, 
one can identify a spin structure $s$ on $M$ with a characteristic
sublink $K$ of $L$ (see section 3.2 for the definition).
 The Reshetikhin-Turaev invariant $\tau(M)$ is defined as a sum over
all colourings (with some coefficients) of the coloured link
invariant of $L$. 
 The refined
Reshetikhin-Turaev invariant $\tau(M,s)$ is defined analogously,
where the sum is taken over
odd colourings  of $K$ and even colourings
 of $L-K$ only. 
It turns out that
$$\tau(M)=\sum_s \tau(M,s)\, .$$

A refinement of the Turaev--Viro invariant $Z(M)$ of a closed
3-manifold $M$ was done
 in two steps. First,  a state sum $Z(M,h)$ for 
 $h\in H^1(M,{\Bbb Z}_2)$ was defined in [TV],
such that
$$Z(M)=\sum_h Z(M,h).$$
Then Roberts [R] constructed an invariant $Z(M,s,h)$ of a closed
oriented
3-manifold $M$ equipped with a spin structure $s$ and 
 $h\in H^1(M,{\Bbb Z}_2)$, such that
$$Z(M,h)=\sum_s Z(M,h,s).$$

As  is well-known (see [Wi], [At]), a
 theory of quantum invariants of closed 3-manifolds is a part
  of 
topological quantum field theory (TQFT), which associates vector
spaces to
closed surfaces and linear
operators  to 3-cobordisms. In this article,
topological quantum
field theories
 extending the quantum invariants of closed 3-manifolds with spin
structure will be referred to as  `spin' TQFT's.

The first spin TQFT was constructed by Blanchet and Masbaum in [BM].
They use an algebraic technique of [BHMV] in order to extend the
 invariants of [Bl], [KM] and [Tu].  Among the results of [BM] are 
the dimension formula for  modules associated to closed connected
surfaces with spin structure and
 the transfer map from the Reshetikhin--Turaev theory to the
spin TQFT.

In this paper we give a different, geometric construction
of a spin TQFT extending  the refined Reshetikhin--Turaev invariants.
Our construction is  parallel to the one given in [T, Chapter 4].
Whence  we briefly recall the construction of Turaev in section 3.1.
We represent  the vector space $V_{(\Sigma_g, \sigma)}$ associated to 
a closed oriented surface $\Sigma_g$ of genus $g$ with spin structure $\sigma$
as a (subspace of a)
vector space generated by `special' colourings of some ribbon
graph $G^g$ (see Fig.1).
 The graph $G^g$ is chosen in such a way that
its regular neighborhood is a handlebody of genus $g$.
`Special' colourings  is  a subset of admissible colourings
of $G^g$, depending on $\sigma$. 
We show that
$$ V_\Sigma = \oplus_\sigma V_{(\Sigma, \,\sigma)},$$
where $V_\Sigma$ is a vector space associated to $\Sigma$ in
the standard  Reshetikhin--Turaev TQFT. 

We define the operator invariant $\tau(M,s)$ of the spin 3-cobordism
$(M,s)$ as follows:
First,  to each connected component $\Sigma_j$ of genus $g_j$
of the boundary of $M$
we glue a
regular neighborhood of the graph $G^{g_j}$,
containing this graph. 
This results in a 
closed 3-manifold $\tilde{M}$ with some ribbon graph, say $G$,
sitting inside. The graph $G$ is a  disjoint union of the
graphs inside  the handlebodies.
Using a surgery  presentation of $\tilde{M}$ along a link $L$, we
show that there is  a one-to-one correspondence between spin
structures on $M$ and characteristic sublinks of $L\cup G$ (see
section 3.2 for the  definition).
Finally, we define  $\tau(M,s)$ as a refined
Reshetikhin-Turaev invariant of the pair  $(\tilde{M}, G)$, where one
sums over odd colourings of the characteristic sublink (determined by $s$)
 and over even colourings of the  other components of $L$. Note that
$\tau(M,s)$ is an element of the vector space generated by the `special'
colourings of $G$.
We  study   gluing properties of  $\tau (M,s)$ and
give an explicit formula for the projector
$$\tau^\sigma: V_\Sigma \to V_{(\Sigma,\, \sigma)}\, .$$
In addition, we
 show,  that for  connected  $\Sigma$,
 the dimension of $V_{(\Sigma,\, \sigma)}$ 
 coincides with the dimension calculated in [BM]. 
The Reshetikhin-Turaev invariant 
of a 3-cobordism $M$
splits as a sum of the refined invariants, i.e.
$$\tau(M)=\oplus_{\sigma} \sum_s \tau(M,s)\, ,$$
where the sum is taken over  $s$ such  that $s|_{\partial M}=\sigma$.
 
\vspace*{0.2cm}

In section 4
  we construct a spin TQFT extending  Roberts'
invariants. In order to do this,
we use
 a modified  state sum operator
$Z(M,G)$ 
of a 3-cobordism $M$  
together with a   3-valent
graph $G$,  which is a subcomplex of a triangulation  of
$\partial M$ (see  [KS], [BD1] and [BD2]). 
This operator
is  equal to the Turaev--Viro
state sum of $M$
with a triangulation of the boundary $\partial M$ given by the graph
dual to $G$.
The advantage  is that $Z(M,G)$ is  a
homotopy invariant of the graph $G$, which can be
effectively   calculated.

 In [BD2] 
 an isomorphism was constructed between the vector space 
$V_{\Sigma_g}$ 
of Turaev--Viro TQFT and the vector space associated to the two copies of
the graph $G^g$. 
Refining this construction,
  we define
the vector space $V_{\Sigma_g}(\sigma, {\rm h})$
associated to a closed oriented surface  $\Sigma_g$
 with spin structure $\sigma$ and  first cohomology
class h, such that
$$V (\Sigma_g)=\oplus_{\sigma,\,{\rm h}}
\;\; V_{\Sigma_g}(\sigma, {\rm h}).$$
Then we construct
the state sum operator
 $Z(M,s,h)$ of a spin 3-cobordism $(M,s)$ with
$h\in H^1(M,{\Bbb Z}_2)$.

Finally, we show  that 
$$V_\Sigma(\sigma, {\rm h})=V_{(\Sigma, \,\sigma)}\otimes
V_{(-\Sigma, \,\sigma+{\rm h})}
$$ and
$$Z(M,s,h)=\tau(M,s)\otimes \tau(-M, s+h),$$
where a negative sign  means the orientation reversal.
This proves that the operator $Z(M,s,h)$ gives rise to  an (anomaly free
non-degenerate) TQFT  on  compact oriented 3-cobordisms equipped
with a spin structure and a first ${\Bbb Z}_2$-cohomology class.


\newtheorem{satz}{Theorem}
\newtheorem{lem}[satz]{Corollary} 
\newtheorem{lemma}[satz]{Lemma}
\input epsf.sty

\section {Initial data and notation}\label{sixj}
In this section
 we define basic algebraic
data, which will be used in the construction 
of invariants.

Let $A$ be a primitive root of unity of order $4r$, where 
$r\in {\Bbb N}$ and $r=0\pmod 4$. 
Consider the  set $I=\{ 0, 1, 2,
...,r-2\}$. For each $i\in I$, we fix  complex numbers
 $\omega_i$ and $q_i$, such that

\be\label{1}
\omega^2_i= (-1)^{i}[i+1]\;\;\;{\rm and} \;\;\;\;q^2_i=(-1)^i A^{i^2+ 2i},\ee
where
$$
 [n]=\frac{A^{2n}-A^{-2n}}{A^2 -A^{-2}},\,\;\;\;\;{\rm for}\;\;
 n\in {\Bbb N}. $$
Furthermore, we choose a complex number $\omega$, such that
\be \label{omega}
\omega^2 =\sum_{i\in I} \omega^4_i= \frac{-2r}{(A^2-A^{-2})^2}\, .
\ee

These data come from  the modular category 
 provided by   `good'
representations of the quantum group
 $U_q({\ma sl}(2,\Bbb C))$ (see  [RT]), where
$A^4=q$. 
In this article we enumerate irreducible
representations of $U_q({\ma sl}(2,\Bbb C))$ 
 by doubled spins $i\in I$. We recall that
 $\omega^2_i$ is equal to the quantum dimension 
of  the $i^{\rm th}$ representation and
the ribbon
graph invariant, defined in [RT],  is multiplied by
$q^{-2}_i$ 
under one  twist on an $i$-coloured  ribbon:
 \begin{center}
\mbox{\epsfysize=2cm
\epsffile{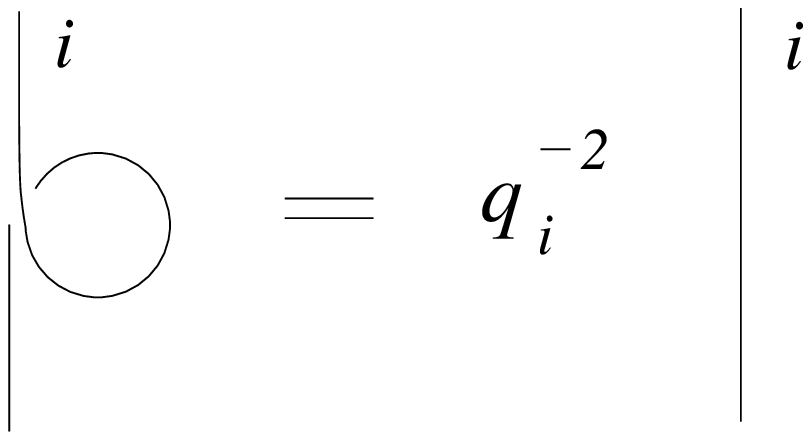}}
\end{center}

A triple $(i,j,k)\in I^3$ is called admissible if
$i+j+k$ is even and 
\be\label{uneq}
 i\leq j+k,\;\,  j\leq i+k,\;\,  k\leq i+j,\;\,  i+j+k\leq 2(r-2)\, .\ee

We finish this section by collecting relations which
 will be of importance in the sequel.
It was shown in [R] that
\be \label{circ1}
\vspace*{1.3cm}
\sum^{r-2}_{i=0,\,i \,{\rm even}} \omega^2_i \hspace*{2cm}
\includegraphics{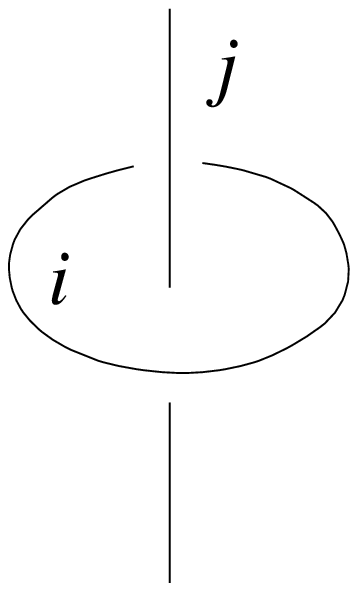}=
\frac{\omega^2}{2}(\delta_{j,\, 0}+\delta_{j,\, r-2})\, ,
\ee
\be \label{circ2}
\vspace*{1.3cm}
\sum^{r-2}_{i=1,\,i\,{\rm odd}} \omega^2_i \hspace*{2cm}
\includegraphics{circ.ps}=
\frac{\omega^2}{2}(\delta_{j,\, 0}-\delta_{j,\, r-2})\, .
\ee
Moreover,
$$ \omega^2=2\sum^{r-2}_{i=0,\, i\;{\rm even}}\omega^4_{i}=
2\sum^{r-2}_{i=1,\,i\, {\rm odd}}\omega^4_{i}.$$
 In addition, we have
\be\label{qpar}
 q^2_{r-2-i}=(-1)^{i+1} q^2_i , \;\;\;\;\;\;
\omega^2_{r-2-i}=\omega^2_i\, .\ee
It follows that 
$$\Delta = \sum_{i\in I} q^2_i \omega^4_i = \sum_{i \,{\rm odd}} q^2_i \omega^4_i\; .$$
Finally,
\be\label{46}
\Delta \,\bar{\Delta}=\omega^2\, ,
\ee
where  $\bar{\Delta}=\sum_{i} q^{-2}_i\omega^4_i$.


\section{ Spin Reshetikhin--Turaev TQFT}
We begin this section by recalling the standard construction of a TQFT
given by  Reshetikhin and Turaev ([RT] and  [T, Chapter 4]). 
After a brief review on spin structures, we discuss a refinement of
this construction determined by a spin structure on a 3-cobordism.

\subsection[Spin Reshetikhin--Turaev TQFT]{Standard model}
Consider a compact oriented 
3-cobordism $M$ with  boundary $\partial M= (- \partial_- M)
\cup \partial_+ M$, where $\partial_- M$ and $\partial_+M$
are the bottom and top bases of $M$, respectively, and
minus means the orientation reversal.
Assume that the boundary of $M$ is
parametrized, i.e., each connected
component  $\Sigma\subset\partial M$ is
supplied with an 
orientation preserving homeomorphism  $\phi:\Sigma_g \to \Sigma\subset\partial_+ M$ or
  $-\phi:-\Sigma_g \to -\Sigma \subset \partial_- M$, 
where  $\Sigma_g$ and  $-\Sigma_g$ 
are the boundaries
 of a standard  oriented
handlebody ($H^+_g, G^g$) and  
an  oppositely oriented   handlebody ($H^-_g, \bar{G}^g$), respectively.

The handlebody ($H^+_g, G^g$) is defined as a regular 
neighborhood in ${\Bbb R}^3$
of the graph $G^g$, depicted in Fig.1, together with the graph 
itself sitting inside.
\begin{center}
\mbox{\epsfysize=2.5cm\epsffile{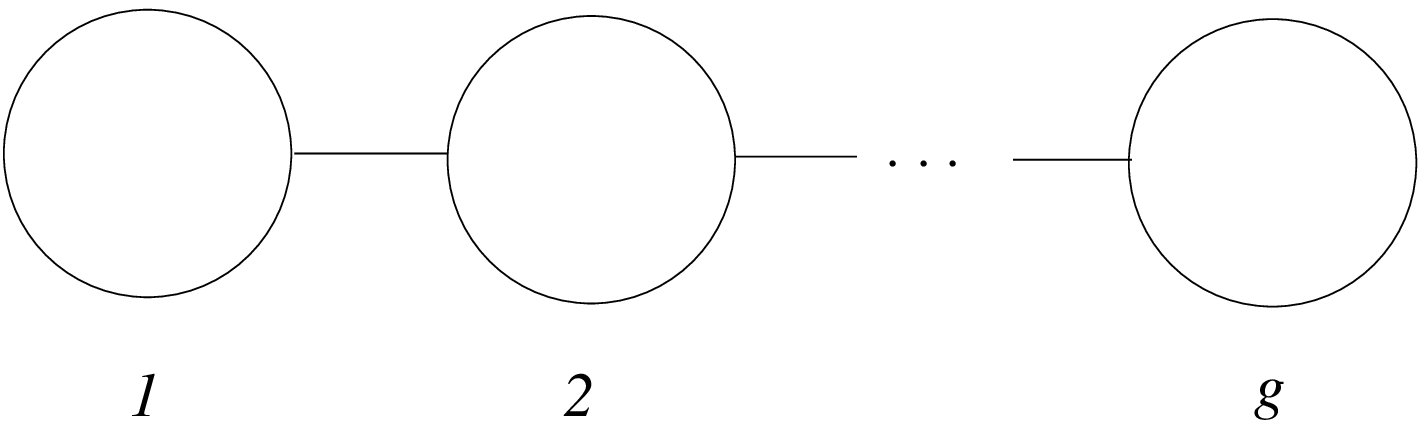}}
\\Fig.1 {\it The 3-valent graph $G^g$ }
\end{center}
The  mirror image of ($H^+_g, G^g$)
 with respect to a horizontal plane in ${\Bbb R}^3$
 defines the  oppositely oriented   handlebody
($H^-_g, {\bar{ G}}^g$).  

By an admissible colouring ${e}=\{e_1, e_2, ..., e_{3g-3}\}$
 of  $G^g$, we mean  an assignment
 of a colour (from $I$) to each line of
 $G^g$, so that  the three colours 
of  lines, meeting in a
3-vertex, form an admissible triple  in the sense of Section 2. We will denote
the $e$-coloured 3-valent graph by $G^g_{{e}}$.
\begin{center}
\mbox{\epsfysize=2.5cm\epsffile{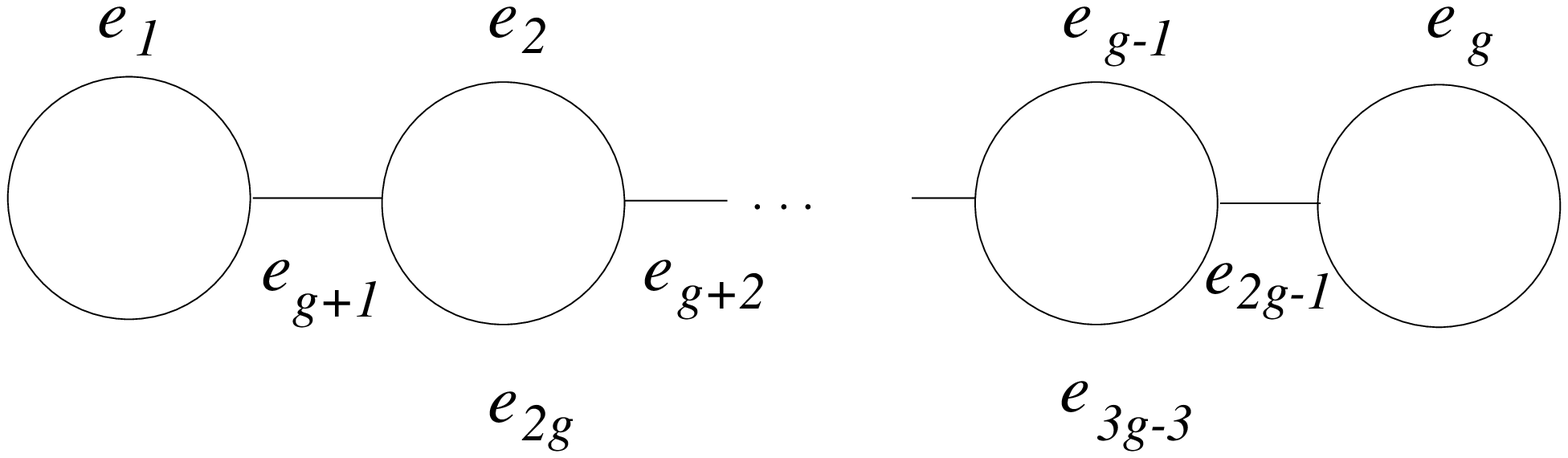}}
\\Fig.2 {\it The  coloured 3-valent graph $G^g_{{e}}$} 
\end{center}
We note that the
admissible colourings of $G^g$  provide  a basis
of the vector space $V_{\Sigma}$ associated by the
Reshetikhin--Turaev TQFT   to a closed parametrized
 surface $\Sigma$ of genus $g$. 
Their number is equal to the  dimension of
$V_{\Sigma}$  given by the  Verlinde formula.
 To a non-connected
surface  one associates  the  tensor product of the vector spaces belonging  to
 connected components. 

\vspace*{0.2cm}

The construction of a 3-cobordism invariant is  as follows:
To each connected component of $\partial_- M$ of genus $g$
one glues a copy of
($H^+_g, { G}^g$) 
along  the given parametrization and analogously one glues
the oppositely oriented  handlebody
to each connected
component of $\partial_+ M$. The result is a closed
3-manifold $\tilde{ M}$ with a ribbon graph, say $G^+\cup G^-$,  sitting
 inside.   The graph $G^+\cup G^-$ is  the disjoint union
of  graphs $\bar{G}^g$ and $G^g$  inside
 the standard handlebodies.
Now the invariant of the 3-cobordism $M$ is defined as an invariant 
of the pair ($\tilde{M}, G^{+}\cup G^-$). 
More precisely, this invariant in the basis, given by the admissible
colourings of $G^+ \cup G^-$, can be written as follows:
\be\label{tau}
\tau(M)_{{e}{e}^\prime}= 
(\Delta \omega^{-1})^{\sigma(L)} \omega^{-m -1+ 
\frac{\chi(\partial_+ M)}{2}} \omega_{{e}\,}
\omega_{{e}^\prime}\sum_{{c}}
\omega^2_{{c}}\;  Z(L_{{c}} \cup G^+_{{
e}}\cup G^-_{{e^\prime}})\, ,\ee
where 
$$\omega_e=\prod_{i} \omega_{e_i},$$
$e$ (resp. $e^\prime$) is a colouring of $G^+$ (resp. $G^-$),
$L\subset S^3$ is an $m$-component surgery link for $\tilde{M}$,
${c}=\{c_1, c_2, ..., c_m\}\in I^m$ 
is a colouring of $L$, $\sigma(L)$ is the  signature
of the linking matrix, $\chi$ is the Euler characteristic
and $Z(G_{{e}} )$
denotes the invariant of a coloured ribbon graph $G_{{e}}$ in $S^3$
as  defined in [RT].

 We set
\be \label{ful}
\tau(M)=\oplus_{{e}{ e}^\prime} \;
\tau(M)_{{e}{e}^\prime}: V_{\partial_- M}\to V_{\partial_+ M}\, .
\ee

It was shown in [T]
that the linear operator $\tau (M)$ determines a TQFT. 
In particular, this means 
that  gluing of  cobordisms  is described by   composing  operators
and that 
$$ \tau(\Sigma\times [0,1])=id_{V_\Sigma}\, .$$

This
 construction 
can be naturally  generalized to 3-cobordisms between
punctured surfaces.
The only significant modification requires the notion of a standard 
handlebody. 

Consider the  handlebody $H^+_g (p)$, whose boundary is an
oriented  surface $\Sigma_g$ with a set ${p}=\{p_1,p_2, ..., p_n\}$
of distinguished points (punctures).
 Attach to each puncture a colour from the 
set ${a}=\{a_1, a_2, ..., a_n\}$ and
embed the graph $G^g({a})$
depicted below
\begin{center}
\mbox{\epsfysize=2.5cm\epsfxsize=10cm\epsffile{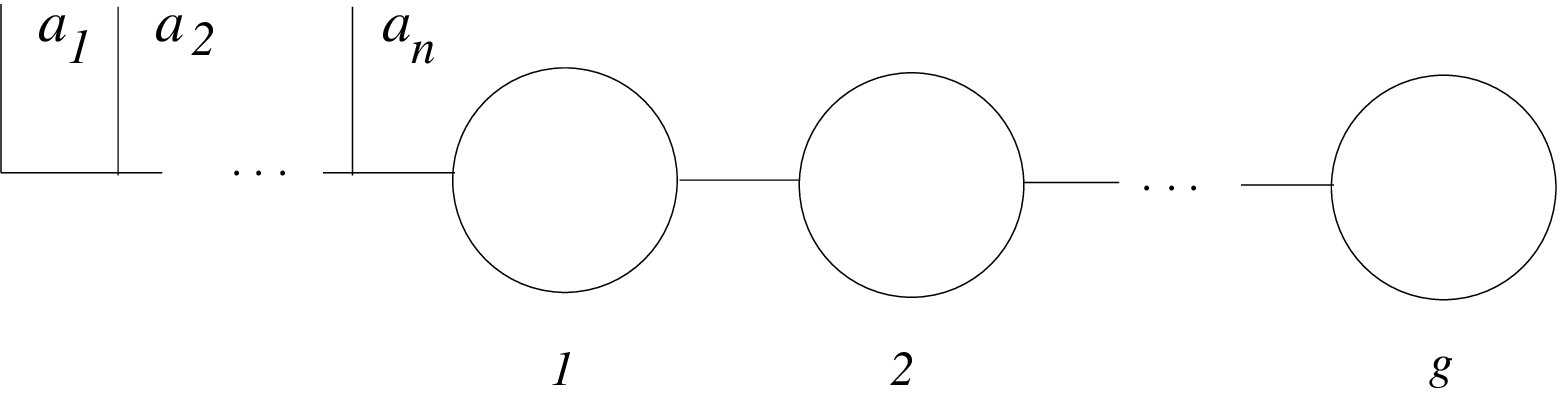}}
\\Fig.3 {\it The  graph $G^g({a})$ }
\end{center}
in $H^+_g (p)$,  so that its  1-vertices 
lie on $\Sigma_g$ and coincide with the punctures $p$ and the remainder of
the graph forms a deformation retract of $H^+_g$.
The resulting pair  $(H^+_g({p}), G^g({a}))$ is a punctured
standard  handlebody. 
A construction of a TQFT is quite analogous  to the one described
above and will  not be repeated here. We mention only that the vector
space 
associated by this TQFT to the
punctured surface $\Sigma_g ({p})$ is generated by colourings of
the graph $G^g({a})$.

\vspace*{0.2cm}

\subsection[Spin Reshetikhin--Turaev TQFT]{Spin structures on manifolds}

A spin structure on an
$n$-dimensional manifold $N$ is a homotopy class of 
a trivialization of the tangent bundle
of $N$ over the 1-skeleton 
which  extends over the 2-skeleton (see
[Ki]). 
The number of different spin structures on $N$
(if it is not zero) is
equal to the number of elements in $H_1(N)$\footnotemark[1].
\footnotetext[1]{Throughout
this paper
all (co)homology groups will have ${\Bbb Z}_2$-coefficients.}
 Moreover,
the whole set of spin structures on $N$ (if it is not empty) is obtained
 by adding  elements of $H^1(N)$ to any fixed spin structure.

There exist two
spin structures on a circle: the bounding spin structure (which extends
over a disc) and the non-bounding or Lie spin structure.   A spin
structure $\sigma$ on a connected
 surface $\Sigma$ defines a
quadratic form $q_\sigma: 
H_1 (\Sigma)\to 
{\Bbb Z}_2$, such that for any embedded closed curve
$\gamma$,  $q_\sigma(\gamma)=0$,
 if $\sigma |_\gamma$ is 
bounding, and $q_\sigma(\gamma)=1$ otherwise (see [Jo]). To
determine a  spin structure on a surface, 
it is  sufficient to say which
simple closed curves in a canonical homology basis (as in Fig.4) are
spin bounding and which are not. 
 
One can also think of 
a spin structure on a manifold $M$  as being
a first
cohomology class of an oriented frame  bundle $F(M)$,
whose restriction to each fibre is non-trivial.  If $M$ is
3-dimensional, this class can be evaluated on a framed knot in $M$,
representing a 1-cycle in $F(M)$
(the rest of a true frame can be reconstructed using
the orientation of $M$). This cohomology class is equal to $1$ on a
trivial knot  in $M$ with zero framing. 

Let us  denote by ${\rm Spin}(M)$
a set of spin structures 
on a 3-manifold $M$. Suppose that $M$ is
obtained by surgery on a framed
 $m$-component  link $L$. Denote by
 $S^3\backslash L$ the 3-sphere
 $S^3$ with a regular neighborhood of $L$ removed.
Then one
can  identify ${\rm Spin}(M)$
with a subset of Spin$(S^3\backslash L)$,
consisting of all  spin structures  which are equal to 1  on 
each component $L_i$ of $L$. 

Taking into account that
$${\rm Spin}(S^3\backslash  L)=s_0 + H^1(S^3\backslash L),$$
where $s_0$ is a spin structure on $S^3\backslash L$,
induced by the unique spin structure on $S^3$,
we observe that any spin structure on $S^3\backslash 
L$ is completely determined  by
its values on the meridians $\{m_i\}^m_{i=1}
$ of the regular neighborhood of $L$. 
One can evaluate a cohomology class $s\in {\rm Spin}(S^3\backslash  L)$
 on a framed knot
$\gamma$ in $S^3\backslash L$ as follows:
$$s(\gamma)= 1+\gamma\cdot\gamma +\sum^m_{j=1} (\gamma\cdot
L_j)(1+s(m_j)), $$
where   $\gamma\cdot L_j={\rm lk}(\gamma,
L_j)$ is the linking number  and
$\gamma \cdot \gamma$ is the framing  on $\gamma$. 
Imposing the condition 
$$s(L_i)=1\;\;\;\;\;\;{\rm  for}\;\;\;\;\;\;\; i=1,2,...,m ,$$ 
we obtain that
any  spin structure $s\in {\rm Spin}(M)$ defines a sublink $K\subset
L$, such that 
 for any component ${L}_i$
of  ${L}$
\be\label{char} {L}_i \cdot {K}={L}_i \cdot {L}_i\, .\ee
The sublink $K$ satisfying (\ref{char}) is called a characteristic
sublink of $L$. 
It  consists of all
the components ${L}_i$ of $L$, such that $s$ is non-bounding on
the  meridian $m_i$ of $L_i$ or, in other words, $s(m_i)=0$.
We define   a characteristic
coefficient $c_i\in {\Bbb Z}_2$  of the component $L_i$ of $L$
equal to one if  $L_i\in K$ and zero
otherwise.

For a 3-cobordism $M$ with parametrized boundary,
 one can identify Spin$(M)$ with a subset of
$${\rm Spin}(S^3\backslash (L \cup
G^+\cup G^-))
=s_0 + H^1(S^3\backslash(L\cup
G^+\cup G^-)),$$
consisting of all spin structures  which are equal to $1$ on $L$ 
 (see Section 3.1 for the definition of $G^+\cup G^-$). 
A basis in $H_1(S^3\backslash
(L\cup G^+\cup G^-))$ is given by 
meridians $\{m_i\}$ of $L$
together with  meridians $\{b_i\}$ of (a
regular neighborhood of)  $G^+\cup G^-$. Denoting by $\{a_i\}$  the longitudes
of $G^+\cup G^-$, we have that
$$s(L_i)=1+L_i\cdot L_i + \sum_j (L_i\cdot L_j)(1+s(m_j))+
\sum_j (L_i\cdot a_j)(1+s(b_j))\, ,$$
where  
$s\in {\rm Spin}(S^3\backslash
(L\cup G^+\cup G^-))$.
 It follows that there exists a one-to-one correspondence between the
solutions of the following equations
$$L_i\cdot (K+A)=L_i\cdot L_i , \;\;\;1\leq i\leq m,$$
where  $K\subset L$ and $A\subset \cup_i a_i$, and 
 spin structures on a 3-cobordism $M$, which do not extend over the
meridians of $K$ and $A$. We will call $K$ a characteristic sublink of
$L\cup G^+\cup G^-$.




\subsection[Spin Reshetikhin--Turaev TQFT]{Spin Reshetikhin--Turaev model}
In this section we construct a spin TQFT
by  refining the model described in section 3.1.

\vspace*{0.2cm}
\centerline{\bf Definition of invariants}
\vspace*{0.2cm}

We start by modifying  the notion
of a  standard handlebody.

Consider the handlebody $H^+_g$ with the boundary $\partial H^+_g=\Sigma_g$
  as depicted below.
 \begin{center}
\mbox{\epsfysize=2.8cm\epsfxsize=8.5cm
\epsffile{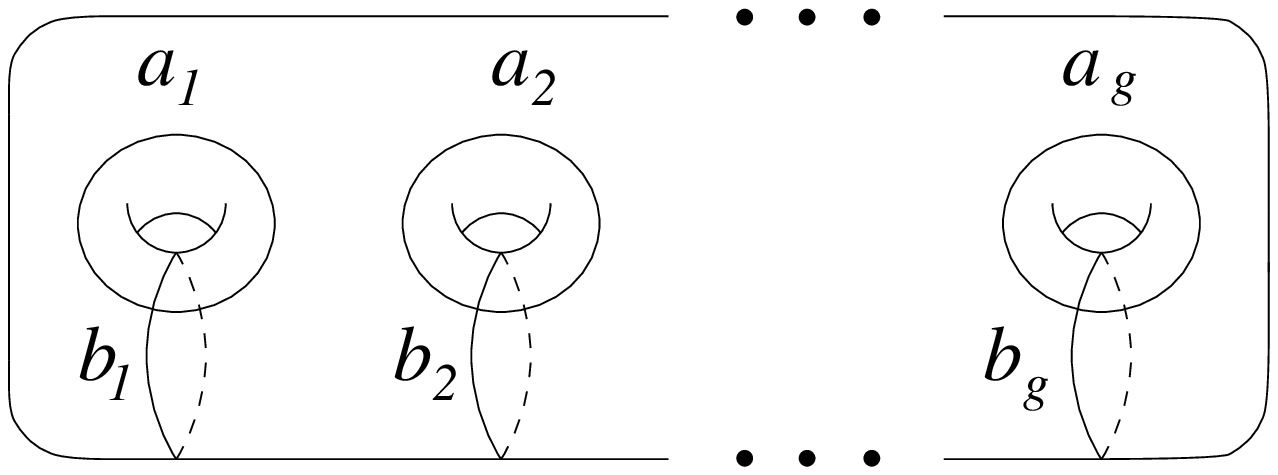}}
\\Fig.4 {\it The canonical homology basis on $\Sigma_g$}
\end{center}
Associate to each meridian $b_i$ of $\Sigma_g$ a number ${\ma s}_i\in
{\Bbb Z}_2$ and denote by $\ma s$ the sequence of these numbers, i.e.
$${\ma s}=\{{\ma s}_1, {\ma
s}_2, ..., {\ma s}_g\} \in {\Bbb Z}^g_2\, .  $$
 Then we embed the graph $G^g$ (see Fig.1)  in $H^+_g$ as its deformation
retract. The resulting triple
($H^+_{g},G^{g}, {\ma s}$) will be called 
a standard handlebody.
 The  oppositely oriented  handlebody 
($H^-_{g}, \bar{G}^{g}, {\ma s}$) 
is defined  by  a mirror image of ($H^+_{g},G^{g}, {\ma s}$). 

Let  $E_{\ma s}$ be 
a subset of  admissible colourings of the graph $G^{g}$
subject to the following relation:
\begin{itemize}
\item
 a colour $e_i\in I$ , $1\leq i\leq g$, is even, if ${\ma s}_i=0$, and
odd otherwise.
\end{itemize}
In the sequel we will  call the elements of
$E_{\ma s}$   {\it special colourings} of the graph
$G^{g}$.
\vspace*{0.2cm}

By a parametrized  surface $(\Sigma, {\ma s})$ of genus $g$  we
understand an oriented closed connected  surface 
of genus $g$ supplied with an orientation
preserving homeomorphism 
$$\phi :\Sigma_g\to \Sigma$$
and a sequence $\ma s$
 of ${\Bbb Z}_2$-numbers associated to $\phi(b_i)$, $1\leq i\leq g$.
We denote by $V_{(\Sigma,{\ma s})}$ the vector space 
associated
to the parametrized  surface $(\Sigma, {\ma s})$,
 which is 
generated by the special colourings
$E_{\ma s}$ of the graph $G^{g}$. Clearly,
$$ V_{\Sigma}=\oplus_{{\ma s}} V_{(\Sigma, {\ma s})},$$
where $V_{\Sigma}$ denotes as before  the vector space associated to
$\Sigma$ in the standard Reshetikhin-Turaev model
and the direct sum is taken over $2^g$ possible choices of $\ma s$. 
To  disjoint
unions of surfaces  we
associate  the tensor product of vector spaces.

\vspace*{0.2cm}
Consider a spin 3-cobordism ($M,s$) with parametrized boundary
$\partial M = (-\partial_- M
)\cup \partial_+ M$, where $s$ is a spin structure on $M$. 
Let us enumerate the connected components of $\partial M$ by an index
$j$, $1\leq j\leq n$. Suppose that the first $l$ of them belong to
$\partial_- M$ and the remaining 
 to $\partial_+ M$. 
Choose a sequence ${\ma s}_j$ of ${\Bbb Z}_2$-numbers associated to
 the $j^{\rm th}$ connected
component $\Sigma_j$ of $\partial M$ in such a way, that
$$({\ma s }_j)_i=q_{s|_{\Sigma_j}}(\phi_j(b_i)),\;\;\;\; 1\leq i\leq g_j,$$
where $\phi_j:\Sigma_{g_j}\to \Sigma_j$
is the  parametrization homeomorphism.

After gluing  (along the parametrizations)
of ($H^+_{g_j},G^{g_j},{\ma s}_j$), $1\leq j\leq l$, 
and  
($H^-_{g_j}, \bar{G}^{g_j},{\ma s}_j $), $l< j\leq n$, 
 to  connected components of $\partial_- M$ and
$\partial_+ M$, respectively,  we
obtain a closed manifold $\tilde M$ with the graph, say 
$G^+ \cup G^-$, sitting
inside.
Denote by $L$ an $m$-component
surgery link for $\tilde{M}$. In general, the spin structure $s$ does
not extend over $\tilde{M}$, but it determines a spin
structure on $S^3\backslash (L\cup G^+\cup G^-)$. Now we
choose a
characteristic sublink $K$ of $L\cup G^+\cup G^-$,
 consisting of all the components $L_i$ of $L$, such that $s$ is 
non-bounding on the corresponding meridians.
Set
$$
\tau(M, s)_{{e}{e}^\prime}= 
(\Delta \omega^{-1})^{\sigma(L)} \omega^{-m -1+ 
\frac{\chi(\partial_+ M)}{2}} \omega_{{e}}\,
\omega_{{e}^\prime}$$
\be\label{stau}
\sum_{{c}\; {\rm odd}}
\omega^2_{{c}} \sum_{{b}\; {\rm even}}
\omega^2_{{b}}\;
 Z(K_{{c}} \cup (L-K)_{{b}}\cup
 G^+_{e}\cup G^-_{{e}^\prime})\; ,\ee
where
$e\in E_{{\ma s}_+}$ 
 and  $e^\prime\in E_{{\ma s}_-}$ are   special colourings of $G^+$ and
$ G^-$,
respectively. Here 
 $${\ma s}_+=\cup^n_{j=l+1} {\ma s}_j , \;\;\;\;  \;\;\;\;
{\ma s}_-=\cup^l_{j=1} {\ma s}_j$$ 
and we  denote by
$c$ and $b$  the colourings of $K$ and $L-K$,
respectively. A colouring is called even (resp. odd), if all its
values are even (resp. odd).

We define the linear operator
$$\tau(M,s): V_{(\partial_- M,\, {\ma s}_-)} \to V_{(\partial_+ M,\, {\ma s}_+)}$$
 corresponding to the spin cobordism
$(M,s)$ by taking a direct sum  over all 
  special colourings of $G^+$ and $ G^-$, i.e.,
\be \label{sful}
\tau(M,s)=\oplus_{{e}{  e}^\prime} 
\tau(M,s)_{{e}{e}^\prime}\, ,\;\;\;\;\;\;
{e}\in E_{{\ma s}_+}, \;\;{e}^\prime\in E_{{\ma s}_-}.\ee

\begin{satz} 
$\tau(M,s)$ is a topological invariant
of a  compact  spin 3-cobordism $(M, s)$ with 
 parametrized boundary.
\end{satz} 

We say that two spin cobordisms $(M,s)$ and $(M^\prime, s^\prime)$ with
parametrized boundary are 
spin homeomorphic if there exists a spin homeomorphism 
$f:(M,s)\to(M^\prime, s^\prime)$  which preserves the parametrized
bases
(or, in other words, whose restriction to the boundary commutes with the
parametrizations). 
\begin{lemma}
Two spin cobordisms $(M,s)$ and $(M^\prime, s^\prime)$ with
parametrized boundary
are spin homeomorphic if and only if
$(L,K,G^+\cup G^-)$ and $(L^\prime, K^\prime, G^{\prime +}
\cup G^{\prime -})$ 
can be
related by (a sequence of) the following refined Kirby move(s).

Add to $L$ an unknotted component $L_i$ 
with framing $\varepsilon= \pm 1$ and characteristic coefficient 
\be\label{9}
c_i=1+
{\rm lk}(L_i,K) +{\rm lk}_{\rm odd} (L_i, G^+\cup G^-)\ee 
and 
change simultaneously
the
part of
$L\cup G^+\cup G^-$, lying in a regular neighborhood of a disc bounded by
$L_i$, 
 by giving a twist (right or left handed, depending on the sign of
$\varepsilon$). 
The last term in (\ref{9}) denotes the linking number
of $L_i$ with the odd coloured lines of the graph $G^+\cup G^-$. 
\end{lemma}

{\bf Proof of Theorem 1:} 
One have to
 show that (\ref{stau}) is invariant under the refined Kirby
move. It is not difficult to verify by direct calculation (see also
 [KM] or [Bl]) that adding of an odd (resp. even)
 coloured unknotted $\varepsilon$-framed component to $L$, linked with
even (resp. odd) number of odd coloured strings 
\footnotemark[2],
\footnotetext[2]{Fusion
preserves the parity of colours.}
and twisting of these strings,
 will multiply the second line in
(\ref{stau}) by $\Delta$ (if $\varepsilon =-1$) or by $\bar{\Delta}$
(if $\varepsilon =1$) 
 and the first line by $\Delta^{-1}$ (if $\varepsilon =-1$) or
by $\omega^{-2}\Delta$ (if $\varepsilon =1$). The claim follows now
from (\ref{46}). 
   $\hfill\Box$

The construction described above can be  straightforwardly
generalized to the case, when the surfaces
$\partial_\pm M$ are provided with punctures 
coloured by $a_\pm$. The corresponding operator invariant is  denoted by
$\tau( a_+, M,s, a_-)$.

\vspace*{0.2cm}
\centerline{\bf
Presentation of spin cobordisms by special ribbon graphs}  
 \vspace*{0.2cm}

In (\ref{stau}) we represented  a  spin 3-cobordism $(M,s)$
by some  special ribbon graph
 $K\cup (L-K)\cup G^+ \cup G^-$ in $S^3$. 
We recall that
 $K$ is the odd coloured,
characteristic  sublink of $L\cup G^+\cup G^-$ and the colourings of 
$G^+$ and
$G^-$
are determined by ${\ma s}_+$ and ${\ma s}_-$, respectively.
It turns out that this construction is 
invertible. This means that each such
special ribbon graph 
 gives rise to a 3-cobordism $M$ with
certain spin structure $s$. Starting from the special ribbon graph, one
can construct ($M,s$) as follows: 


One removes  tubular neighborhoods
of $G^+$ and $G^- $ from $S^3$. This results in a 3-cobordism $E$ with
bottom base 
 $\Sigma^-$ and top base
 $\Sigma^+$. We  provide $E$  with 
orientation induced by right-handed orientation in $S^3$ and bases
with orientations, such that $\partial E= (-\Sigma^-)\cup \Sigma^+$.
 We choose the parametrizations, 
which send the $a$-cycles of $\Sigma_g$
to the loops  on $\Sigma^\pm $ 
homotopic to the circles of the graphs $G^\pm$.
  Now remove from $E$ a regular neighborhood of $L$.
Choose a spin structure $s$ on $E\backslash L$, which is non-bounding on the
meridians of $K$ and on the meridians  of the
odd coloured lines of $G^{+}\cup G^-$. Glue solid
tori back to $E\backslash L$ 
along the homeomorphisms determined by framing.
This results in an oriented 
3-cobordism, say $M$, with spin structure $s$ and parametrized boundary.

\vspace*{0.2cm}

\centerline{\bf Gluing properties}
\vspace*{0.2cm}
We will show that 
the operator $\tau (M,s)$ defines  a non-degenerate spin TQFT. 

\begin{satz}
If the spin 3-cobordism $(M,s)$ is obtained from  $(M_1, s_1)$
and $(M_2, s_2)$  by gluing
along a homeomorphism $f :\Sigma\to\Sigma^\prime$
which preserves
spin structures and
commutes with parametrizations, then 
\be \label{glu}
\tau(M,s)_{{e} {e}^\prime}=
k\sum_{{e}^{\prime\prime}\in E_{\ma s}}
\; \tau(M_2,s_2)_{{e} {e}^{\prime\prime}}\;
\tau(M_1, s_1)_{{e}^{\prime\prime}
 {e}^\prime}\, ,\ee
where $\Sigma= \partial_+ M_1$,
$\Sigma^\prime = {\partial_- M_2}$ are parametrized connected surfaces and
$k= (\Delta \omega^{-1})^{\sigma(L)-\sigma(L_1)-\sigma(L_2)}$ is
an anomaly factor.  
 \end{satz}

{\bf Proof:} We can
represent $M_1$ and $M_2$ by special ribbon graphs
 $K_1\cup (L_1-K_1)\cup \bar{G}^g \cup G^-_1$
and  $K_2\cup (L_2-K_2)\cup G^+_2 \cup G^g $, 
respectively, where $g$ is the genus of $\Sigma$. 
 Putting the special ribbon
graph representing $M_2$ on the top of the
graph for $M_1$ and
summing over $e^{\prime \prime}_i$ ($e^{\prime\prime}\in
E_{\ma s}$)
with $i>g$,
 we obtain a
special ribbon graph 
\be \label{gr} K_2\cup K_1\cup (L_2-K_2)\cup( L_1-K_1)\cup \Omega
\cup G^+_2 \cup G^-_1 
 \, ,\ee
where by $\Omega$ we denote the
$g$ annuli, which remain of
$ G^g$ and $\bar{G}^g$ after the  summation. 
The  graph (\ref{gr}) is, in fact, a special ribbon 
graph representing $M$ (see [T,
p.175] for more details). Its characteristic sublink consists of
$K_1\cup K_2$ together with the odd coloured annuli of
$\Omega$.
$\hfill\Box$
\vspace*{0.2cm}

{\bf Remark:} Theorem 3 can be straightforwardly  generalized to the
case, when  $\Sigma\subset \partial_+ M_1$,
$\Sigma^\prime \subset {\partial_- M_2}$.
\vspace*{0.2cm}

If we glue  3-cobordisms along non-connected surfaces, 
the situation becomes
more complicated,  because  a spin structure on the resulting
manifold is not uniquely determined by the
spin structure on  3-cobordisms
glued together. In this case we have the following theorem:
\begin{satz}
If the spin 3-cobordism $(M,s)$ is obtained from  $(M_1, s_1)$
and $(M_2, s_2)$  by gluing along a homeomorphism $f :\partial_+ M_1 
\to \partial_- M_2 $  which preserves spin structures  and 
commutes with parametrizations, then 
\be \label{glu1}
\sum_s \tau(M,s)_{{e} {e}^\prime}=
k\sum_{{e}^{\prime\prime}}
\; \tau(M_2,s_2)_{{e} {e}^{\prime\prime}}
\tau(M_1, s_1)_{{e}^{\prime\prime}
 {e}^\prime}\, ,\ee
where the  sum on the left hand side
 is taken over spin structures  such that 
$s|_{M_1}=s_1$ and $s|_{M_2}=s_2$.
 \end{satz}    

{\bf Proof:} Assume that $\partial_+ M_1$  consists of $n$ connected
components of genera $g_1$, $g_2$, ... and $g_n$. Now the special ribbon
graph representing $M$ can be obtained from (\ref{gr}) by replacing
$\Omega$ with a family of $\Omega_i$, $1\leq i\leq n$, 
where by  $\Omega_i$ we denote  $g_i$ annuli,
 and then by encircling   $\Omega_i$, $1\leq i\leq n-1$,  by an
unknotted annulus (see [T, p.177] for more details).
Using fusion rules,
(\ref{circ1}) and (\ref{circ2}) one can split this graph for
$M$ into two parts. The first one consists of a
disjoint union of the special ribbon graphs 
representing
$M_1$ and $M_2$. The second part 
contains terms where
 the special ribbon graph for $M_1$ and $M_2$ are 
 connected by 
$(r-2)$-coloured lines.
 The sign of these terms  depends on
the choice of a spin structure $s$ on $M$, whose restrictions to 
 $M_1$
and  $M_2$ are equal to $s_1$ and $s_2$, respectively. 
Taking the  sum over all $2^{n-1}$  such $s$, we obtain
(\ref{glu1}).
$\hfill\Box$
\vspace*{0.2cm}

In
the next lemma we calculate the invariant of a spin 3-manifold
obtained from two other spin  manifolds by gluing  along a
non-connected surface.

\begin{lemma}
Let $(M,s_i)$,  $i\in {\Bbb Z}_2$, be  spin 3-cobordisms
 obtained from  $(M_1, s_1)$
and $(M_2, s_2)$  by gluing along a homeomorphism $f :\partial_+ M_1 
\to \partial_- M_2 $  which preserves spin structures and 
commutes with parametrizations. Here $\partial_+
M=\Sigma_1\cup\Sigma_2 $, $s_i|_{M_1}=s_1$,  $s_i|_{M_2}=s_2$, $s_0$ 
is bounding and $s_1$ is not bounding 
 on the additional cycle, which appears after
gluing along a non-connected surface. Then 
$$
 \tau(M,s_i)_{{e} {e}^\prime}=
k/2\,\,\,\bigg[ \,\, \sum_{{e}^{\prime\prime}}
\; \tau(M_2,s_2)_{{e} {e}^{\prime\prime}}
\tau(M_1, s_1)_{{e}^{\prime\prime}
 {e}^\prime}\,\,\, + $$ 
\be\label{pun1}
 +\,\, (-1)^i 
\sum_{{e}^{\prime\prime}}
\; \tau(M_2,s_2, r-2, r-2)_{{e} {e}^{\prime\prime}}\;\;
\tau(r-2,r-2, M_1, s_1)_{{e}^{\prime\prime}
 {e}^\prime}\,\bigg],
\ee
where,
 in the second term, 
 $\Sigma_1$ and $\Sigma_2$ are supposed to have an
(r-2)-coloured puncture.  
 \end{lemma}    

{\bf Proof:}  As 
explained in the proof of theorem 4, the special ribbon graph
representing $M$ looks as  follows:
 \begin{center}
\mbox{\epsfysize=2.8cm\epsffile{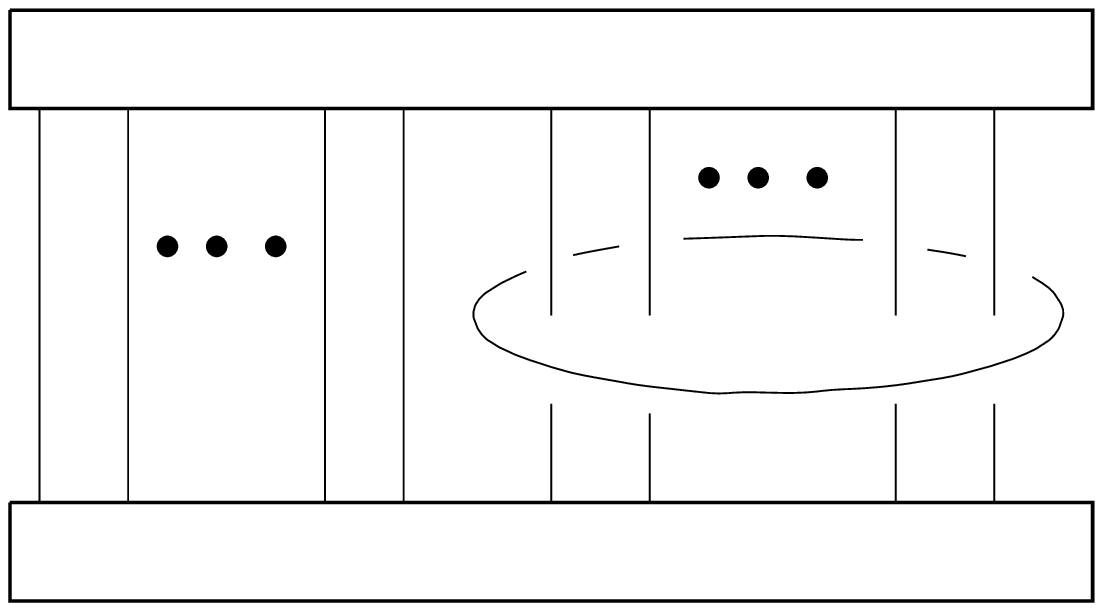}}
\end{center}
where the 
rectangles designate
 the remainder of the ribbon graph. The
circle is odd coloured for $s_1$ and even for $s_0$. Using fusion
rules, one can change this graph in the following way:
 \begin{center}
\mbox{\epsfysize=3.4cm\epsffile{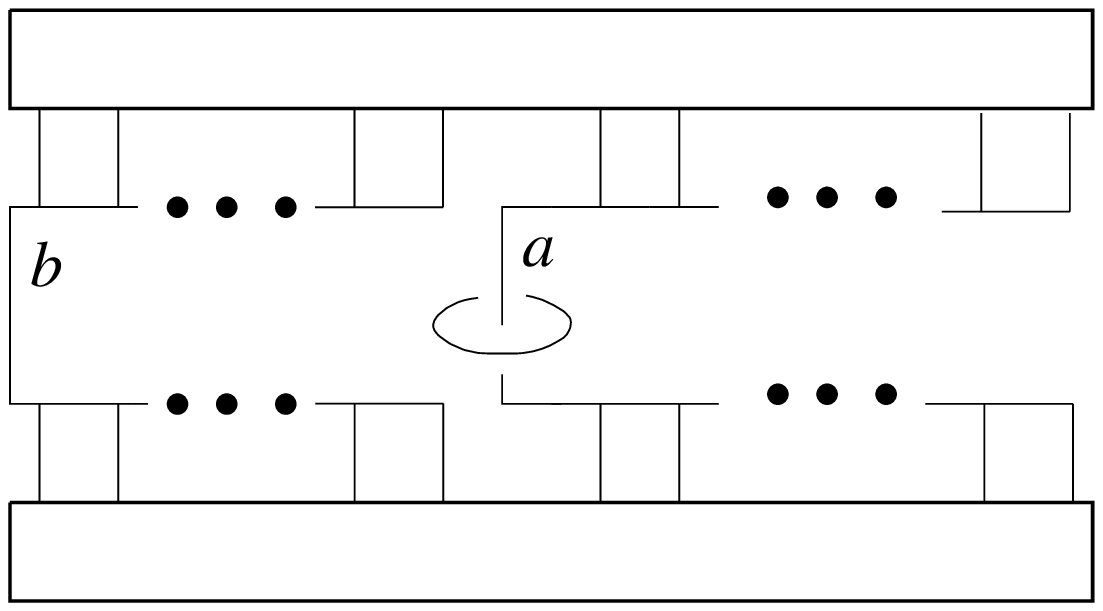}}
\end{center}
where one takes a sum over colourings of the new lines
with quantum dimensions as coefficients. It follows
from (2.4) or (2.5) that the colour $a$ could be either $0$ or $r-2$.
If $a=0$ (resp. $a=r-2$), $b$ should be  equal to $0$ (resp. $r-2$)
too, and we get the first (resp. second) term in (\ref{pun1}).
$\hfill\Box$
\vspace*{0.2cm}

\centerline{\bf Vector spaces associated to surfaces with spin structure}
\vspace*{0.2cm}

 Due to Theorem 3, 
for a spin 3-cobordism $(M,s)$ whose boundary $\partial M=\Sigma$ is a
parametrized
   surface of genus $g$ with spin structure
$\sigma=s|_{\Sigma}$, 
\be \label{56}
\tau(M,s)_{{e}}=\sum_{{e}^\prime\in E_{\ma s}}
\tau(\Sigma\times [0,1],
\sigma\cup\sigma)_{{e}{
e}^\prime}\; \tau(M,s)_{{e}^\prime}\,.\ee


One can now  define the vector space
$V_{(\Sigma, \, \sigma)}$,  associated by the spin TQFT
to the  surface $\Sigma$ with  spin structure $\sigma$, as 
the support of the projector
$$\tau(\Sigma\times [0,1], \sigma\cup\sigma):
V_{(\Sigma,
{\ma s})} \to V_{(\Sigma , {\ma s})}.$$

Assume (without loss of generality) that the parametrization of
$\Sigma$ in the cylinder
 is given by the identity homeomorphism. Then
the special ribbon graph in Fig.5  
 \begin{center}
\mbox{\epsfysize=2.5cm\epsffile{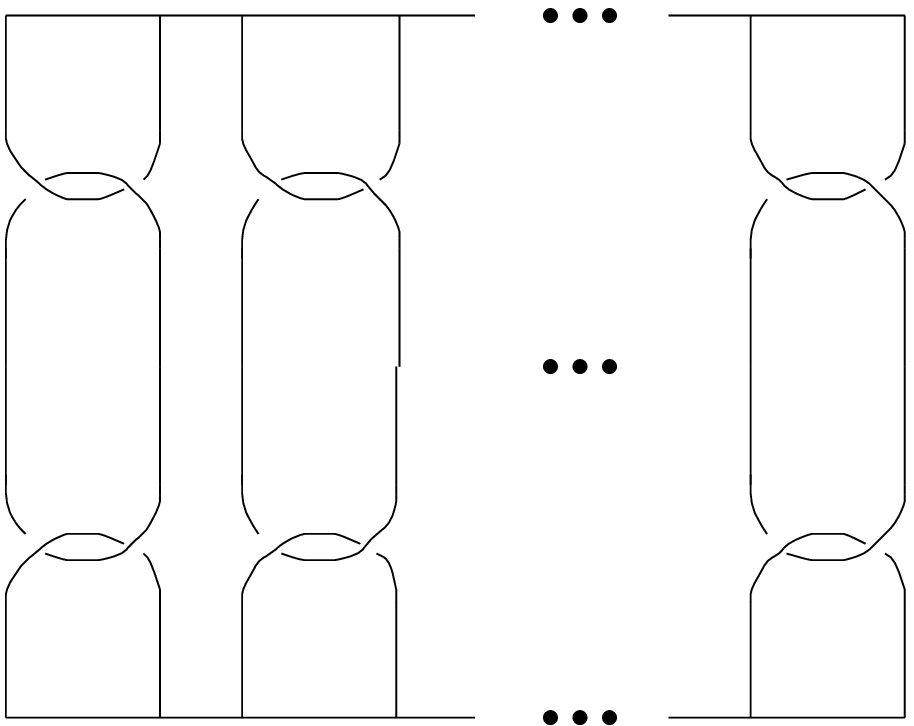}}
\\Fig.5 {\it The special ribbon graph corresponding to a cylinder}
\end{center}
represents the 3-cobordism $\Sigma\times [0,1]$
(see  [T, p.173] for more details). It consists of two
copies of the graph $G^g$ linked with $g$ annuli 
$\Gamma_1\, ,\;\Gamma_2\, ,\; ...\, , \;\Gamma_g$. A  special
colouring of $G^g$ is determined by $q_\sigma(b_i)$,
$1\leq i\leq g$,  and the parity of
colours on $\Gamma_i$ by $q_\sigma (a_i)$.
Using (\ref{circ1}), (\ref{circ2}) and fusion rules  one can
calculate that
\be \label{cyl}
\tau(\Sigma\times [0,1], \sigma\cup\sigma)_{{
e}{e}^\prime}= \frac{1}{2^g} 
(\delta_{e_1\, e^\prime_1}+(-1)^{c_1} \delta_{e_1\, \hat{e}^\prime_1})
...
(\delta_{e_g\, e^\prime_g}+(-1)^{c_g} \delta_{e_g\, \hat{e}^\prime_g})
\prod_{i>g}\delta_{e_i \, e^\prime_i}\, , \ee
where $c_i=q_\sigma(a_i)$ and $\hat{e}_i=r-2-e_i$. 
For simplicity
 we will write $\tau^\sigma$  for
$\tau(\Sigma\times [0,1], \sigma\cup\sigma)$ 
in what follows.
 
One can easily establish that the operators
$\tau^\sigma$ 
form  a family of $4^g$ orthogonal projectors on the vector spaces 
$V_{(\Sigma, \,\sigma)}$, i.e.
$$\sum_{{e}^\prime}
\tau^{\sigma_1}_{{e}{e}^\prime}
\; \tau^{\sigma_2}_{{e}^\prime {
e}^{\prime\prime}}=\cases{0,& if $\sigma_1\neq \sigma_2$\cr
\tau^{\sigma_1}_{{e} {
e}^{\prime\prime}}, & if $\sigma_1=\sigma_2$ \cr}
$$
and 
\be\label{14}V_{\Sigma}=\oplus^{4^g}_{i=1} V_{(\Sigma,\, \sigma_i)}.\ee
As usual, we associate the tensor
product of vector spaces to the  disjoint union of  surfaces. 

Clearly,
$$ \tau(M,s): V_{(\partial_-M,\, s_-)}\to  V_{(\partial_+ M,\, s_+)},$$
where  $s_\pm=s|_{\partial_\pm M}$.

\begin{lemma}
The Reshetikhin--Turaev
 invariant of a 3-cobordism  $M$ with parametrized boundary
splits
as  a sum of the refined invariants corresponding to different spin
structures, i.e. 
\be \label{09}\tau(M)=\oplus_{s_\pm}
\sum_{s}\tau(M,s) \, ,\ee
where the  sum is over all spin structures $s$ on $M$  such that
$s|_{\partial_\pm M}=s_\pm$.
\end{lemma}

{\bf Proof:} 
The claim follows from the fact that the contribution
to $\tau(M)$ from  odd coloured, non-characteristic sublinks 
of $L\cup G^+\cup G^-$  is equal to zero.  The explicit computations
are quite analogous to the one given in  [Bl] or [KM],  and they will
not be repeated here.
$\hfill\Box$

\vspace*{0.2cm}

\centerline{\bf Dimension of vector spaces}
\vspace*{0.2cm}

The Reshetikhin-Turaev
 TQFT  yields a representation of the mapping class group (MCG).
 The matrix elements for  generators of MCG 
 are listed, e.g., in [KSV]. In the spin TQFT, the MCG generates
  transformations between
vector spaces associated to  different spin structures with the
same Arf-invariant. We recall that the Arf-invariant of a quadratic form
$q_\sigma$ (corresponding to spin structure $\sigma$ on $\Sigma_g$) is
defined as follows:
$${\rm Arf}(\sigma)= \sum^g_{i=1}
q_\sigma(a_i)q_\sigma(b_i)\, ,$$
where  $a_i, b_i$ is the symplectic homology basis  depicted
in Fig.4.


As a result, the dimension of $V_{(\Sigma, \sigma)}$ depends only on the
Arf-invariant of $\sigma$, but not on $\sigma$ itself. 
 On $\Sigma$ there exist $2^{g-1}(2^g+1)$ spin structures with
 Arf-invariant equal to zero and $2^{g-1}(2^g-1)$ with
 Arf-invariant equal to one. 

\begin{satz} For a closed  surface $\Sigma$
of genus $g$ with spin structure $\sigma$,
\be\label{11} {\rm dim} V_{(\Sigma,\, \sigma)}=\frac{1}{4^g}\; [\;{\rm dim}
V_{\Sigma} + ({r/2})^{g-1} (2^g-1)\;], \;\;{\rm if
\;\;Arf}(\sigma)=0\,  , \ee
\be\label{2} {\rm dim} V_{(\Sigma, \,\sigma)}=\frac{1}{4^g} [\;{\rm dim}
V_{\Sigma} - (r/2)^{g-1} (2^g +1)\;], \;\;{\rm if
\;\;Arf}(\sigma)=1\,,\ee
where ${\rm dim} V_\Sigma$ is given by the Verlinde formula.

\end{satz}

The dimensions of  spin modules were first calculated in [BHMV]
using a rather developed algebraic technique. Here
we will use   simple
geometric  arguments, which  refine  Lickorish's
calculations in  [Li].
\vspace*{0.2cm}

{\bf Proof:}
The dimension of the vector space $V_{(\Sigma,\, \sigma)}$
can be calculated as follows
\be \label{dim}
{\rm dim} \, V_{(\Sigma,\,\sigma)}= 
{\rm tr}\; \tau (\Sigma\times [0,1],\sigma\cup\sigma).\ee
Theorem 4 implies that
\be \label{6}
{\rm tr}\; \tau (\Sigma\times [0,1],\sigma\cup\sigma)=
\tau(S^1\times \Sigma, s_0) +\tau(S^1\times \Sigma, s_1)\, ,\ee
where $s_i|_{\Sigma}=\sigma$,
 $s_0$  is   bounding and $s_1$ is not bounding  on $S^1$.
A surgery diagram for $S^1\times \Sigma$ can be obtained by taking
$g$ copies of the annulus  containing  a link, which is depicted below,
\be \label{ann}
\mbox{\epsfysize=2.7cm\epsffile{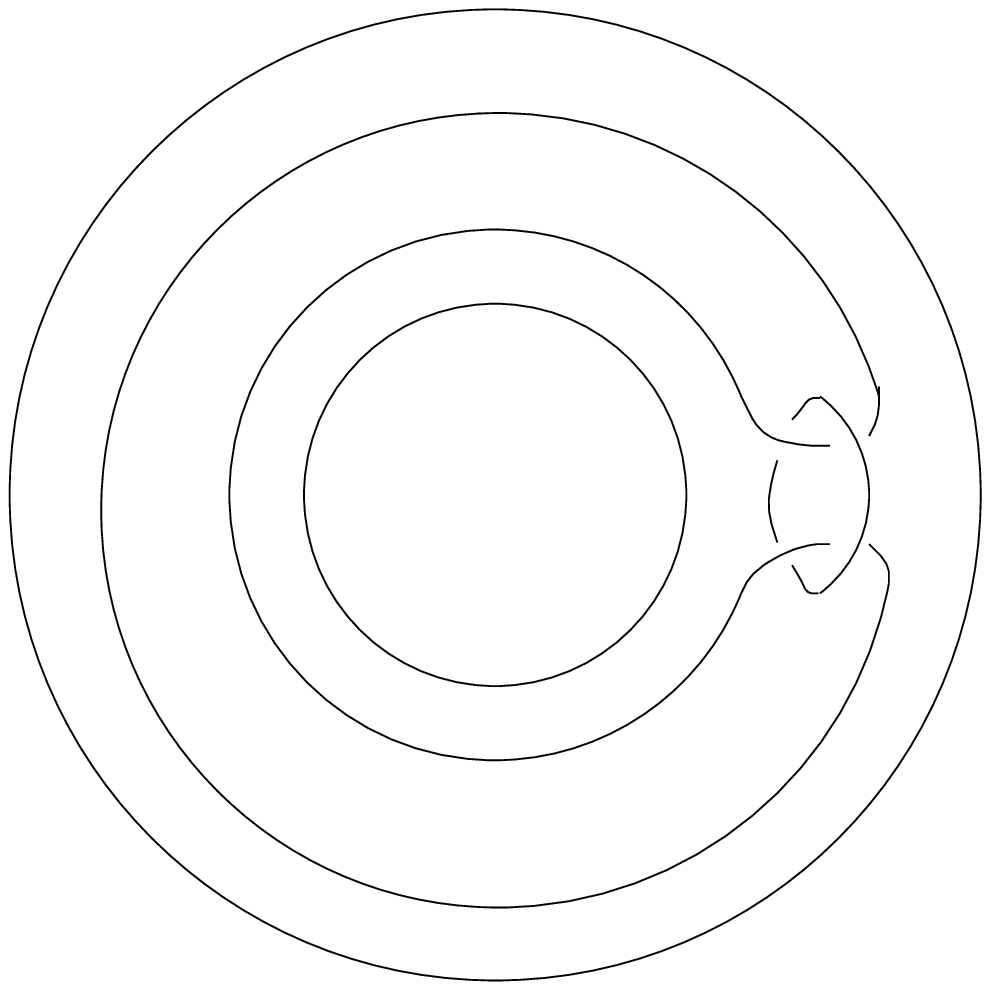}}
\ee
threading un
unknotted closed curve $l$ though the annuli and finally
taking the resultant link
of $2g+1$ components. 

Denote a colour of $l$ by $a$.
Then the invariant $\tau(S^1\times \Sigma, s_i)$
can be calculated in the following way:
One takes $g$ times  expression (\ref{exp}), 
\be \label{exp}
\mbox{\epsfysize=2.5cm\epsffile{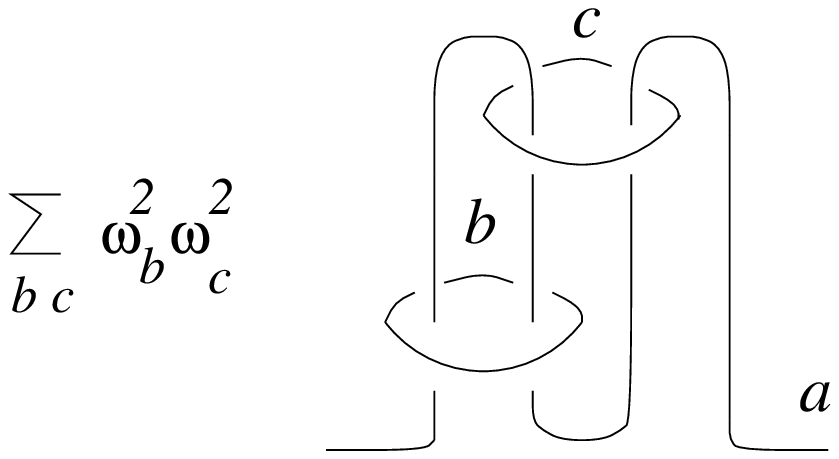}}
\ee
closes an $a$-coloured line,
 sums over $a$ with $\omega^2_a$ as coefficients,
(note that $a$ is even for $s_0$ and odd for $s_1$) and multiplies
by $\omega^{-2g-2}$. 

Consider at first the case  when ${\rm Arf}(\sigma)=0$.
Then one can suppose that all colours (except of $a$) are even.  
Applying fusion rules,  (\ref{circ1}) and
the following  formula
$$\left|\begin{array}{ccc} r/2-1&r-2&r/2-1\\
r/2-1&r-2&r/2-1\end{array}\right| = -\;\omega^{-2}_{r/2-1}$$
(see (4.5) in [BD1] for the graphic  and [TV] for 
the analytic definition of 6j-symbols),
one can reduce
(\ref{exp}) to the $a$-coloured line multiplied by
$$\frac{\omega^4}{ 4\,\omega^4_a}(1+ 
\delta_{a,\,  {r/2-1}})\,\includegraphics{a.ps}
 ,$$ 
where the last term contributes  to $\tau(S^1\times \Sigma, s_1)$ only,
because  $r/2-1$ is odd. Taking into
account all coefficients,
 we obtain that
$$\tau(S^1\times \Sigma, s_0)=\frac{\omega^{2g-2}}{4^g}\sum_{a\, {\rm
even}} \omega^{4-4g}_a\, ,$$
$$\tau(S^1\times \Sigma, s_1)=\frac{\omega^{2g-2}}{4^g} (\sum_{a\, {\rm
odd}} \omega^{4-4g}_a\, +\frac{2^g-1}{\omega^{4g-4}_{r/2-1}})
,$$
which after substituting in (\ref{6}) and using (\ref{1}) and (\ref{omega})
implies (\ref{11}).

The dimension of $V_{(\Sigma, \sigma)}$ with ${\rm
Arf}(\sigma)=1 $ can be calculated analogously or determined from the
formula:
$${\rm dim}\, V_\Sigma = 2^{g-1}(2^g +1){\rm dim} \, V_{(\Sigma, \,
\sigma_0)}+ 
2^{g-1}(2^g -1){\rm dim} \, V_{(\Sigma, \, \sigma_1)}\, ,$$
where  ${\rm Arf}(\sigma_0)=0$ and  ${\rm Arf}(\sigma_1)=1$.
$\hfill\Box$ 

\section{Refined Turaev--Viro TQFT}
The aim of this section is to refine the construction of Turaev--Viro
3-cobordism invariants as given  in [BD1], [BD2]  and define the state sum
operator  $Z(M,s,h)$, satisfying the requirements of a TQFT,
where $s$ is a spin structure on 
$M$ and $h\in H^1(M)$.
We start by recalling the construction of [BD1,2].

\subsection[Refined Turaev--Viro TQFT]{Standard  model}

The Turaev-Viro state sum is defined for any compact triangulated
3-manifold $M$ as follows: One puts colours on 1-simplexes of  $M$ and 
associates 6j-symbols to coloured tetrahedra. 
Then the Turaev-Viro invariant is given by a sum 
over all colourings
of 1-simplexes in the interior of $M$
 of the product of 6j-symbols (with some coefficients).
 The vector
space $V(\Sigma)$
associated to a triangulated
surface $\Sigma$
  is defined as a direct sum over
all colourings  of  the 
tensor product 
 of vector
spaces belonging to  coloured
triangles  of $\Sigma$
modulo some equivalence relation.

As was already mentioned in the introduction,  
we will use
 a modified state sum operator $Z(M,G)$, where $G$ is a
3-valent ribbon graph on $\partial M$. The operator $Z(M,G)$ was
defined in [BD1] (see also [KS]) in such a way, that it is equal to
the Turaev-Viro state sum for $M$,
where the triangulation of  $\partial M$ is
given by the graph dual to  $G$ \footnotemark[3].
\footnotetext[3]{In this article we suppose   that
the graph  $G$ is  large enough
in order that its dual defines  a triangulation of $\partial M$.}
Moreover, $Z(M,G)$ is a homotopy invariant of the graph $G$.
In [BD2] an isomorphism was constructed
between $V(\Sigma)$ and the vector space generated  by
colourings of  two copies of the graph, depicted in Fig.1.

The cobordism $M^+_g$ providing this isomorphism we will call a
standard handlebody. $M^+_g$ is a cylinder
$\Sigma_g\times [0,1]$, where $\Sigma_g$ is a closed oriented surface
of genus $g$
standardly embedded in ${\Bbb R}^3$. Furthermore, $M^+_g$ contains
 an arbitrary $3$-valent  
 graph 
${\cal G}^g$, 
sitting
on $\Sigma_g = \Sigma_g\times \{1\}$,
 and the coloured graph
$G^g_{e}\cup \bar{G}^g_{f} \cup m_{{x}}$,
depicted below, on $-\Sigma_g=\Sigma_g\times \{0\}$,
 \be \label{basis}
\centerline{\hbox{
\psfig{figure=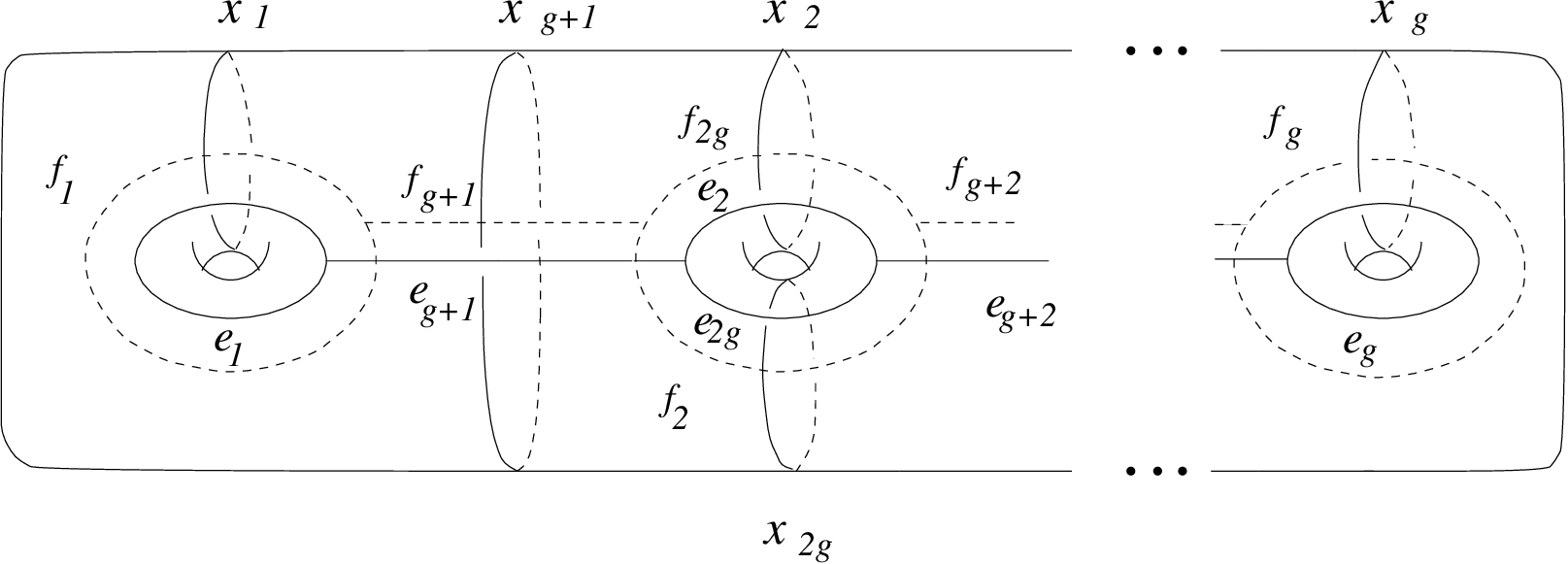,height=6.5cm,width=14.4cm}}}
\ee
where  $m=\{m_1, ..., m_{3g-3}\}$ is the ordered set of meridians
coloured by  $x=\{x_1, ..., x_{3g-3}\}$ and
$e=\{e_1, ...,e_{3g-3}\}$, $f=\{f_1,..., f_{3g-3}\}$ are admissible
colourings of $G^g$
and $\bar{G}^g$, respectively.  We note that
 the ${f}$-coloured graph is drawn on the backward side of $\Sigma_g$.

  The 
 state sum $K_{{e}{f}}$ of the standard handlebody  is
given by the formula:
\be K_{{e} {f}} =
\omega^{g-1}\omega_{{e}}
\omega_{{f}}
\sum_{{x}}
\prod^{3g-3}_{i=1} \frac
{\omega^2_{x_i}}{\omega^2}\; Z(M^+_g,
G^g_{e}\cup \bar{G}^g_{f} \cup m_{{x}} \cup {\cal G}^g),\ee
where  the  sum is over  colourings of the meridians. This state sum
defines a linear
operator 
$$ K_{{e}{f}}: V^L_g({e}) \otimes
V^R_g({f}) \rightarrow V({\Sigma_g}) , $$
where $  V^L_g({e}) \otimes
V^R_g({f})$ is the vector space associated to the graph
$G^g_{e}\cup \bar{G}^g_{f} $.
 It turns out, 
that  the  mirror image $M^-_g$
 of $M^+_g$
yields
an inverse operator
$$ L_{ 
{e}{f}}: V({\Sigma_g}) \rightarrow
V^L_g({e})\otimes V^R_g({f}),\;  $$
 which  satisfies the following equation
 (see [BD2] for more details):
$$L_{{e}^\prime{f}^\prime}
K_{{e}{f}} 
=
\delta_{{e},{e}^\prime}\delta_{{f},{f}^\prime}\;
{\eins}_{V^L_g({e})\otimes V^R_g({f})}.$$
Taking into account that the dimensions of $\oplus_{{e}{f}}\{
V^L_g({e}) \otimes V^R_g({f})\}$ and $
V({\Sigma_g})$ coincide, we obtain that 
$$K=\oplus_{{e}{f}}K_{{e}{f}}
$$
is an isomorphism and
admissible colourings of
$G^g_{e}\cup \bar{G}^g_{f}$ provide a  basis of
$V({\Sigma_g})$. 

From now on
 we fix  the standard handlebodies $M^+_g$
and $M^-_g$ together with the graphs on their boundaries.
We say that an
oriented  triangulated surface $\Sigma$ is parametrized, if
it is supplied with a simplicial map $\phi: ({\cal G}^g)^\ast \to X$,
where by $({\cal G}^g)^\ast$ we denote the triangulation of
$\Sigma_g$, given by the graph dual to ${\cal G}^g$, and $X$ is a
triangulation of $\Sigma$. The parametrization of $-\Sigma$ is given
by the map $-\phi:(\bar{\cal G}^g)^\ast \to -X$.  
\vspace*{0.2cm}

Consider a 3-cobordism $M$ with parametrized boundary
 $\partial M=(-\partial_- M) \cup\partial_+ M$. 
Let us glue the standard handlebodies 
to the connected
components of $\partial_\pm M$  along the
parametrizations.
The state sum of the resulting manifold  with a 3-valent graph on the
boundary defines an invariant of the 3-cobordism $M$ in the 
basis mentioned above. More precisely,
\be \label{Z}
Z(M)_{{e}{f}, {e}^\prime {f}^\prime
}=\omega^{\frac{-\chi(\partial
M)}{2}} \omega_{{e}} \omega_{{f}}
\omega_{{e}^\prime} \omega_{{f}^\prime}
\sum_{{x} {y}} \prod_{i\, j} \frac{\omega^2_{x_i}
\, \omega^2_{y_j}}{\omega^2\, \omega^2} \; 
Z(M, G^+_{{e}{f}} \cup
G^-_{{e}^\prime{f}^\prime}\cup
m_{{x}}\cup m_{{y}})\, ,\ee
where $G^+_{{e}{f}}=G^+_e\cup \bar{G}^+_f$
and $G^-_{{e}{f}}= G^-_e\cup \bar{G}^-_f$ 
are the  disjoint unions of the graphs $\bar{G}^g_{e}\cup
G^g_{f}$ and
$G^g_{e}\cup \bar{ G}^g_{f}$, sitting
on the boundaries of the standard handlebodies $M^-_g$ and $M^+_g$,
respectively. 
Representing $M$ by surgery on an $m$-component link $L$ and using the
technique developed in [BD1] and [BD2], one can rewrite (\ref{Z}) in
terms of the link invariants:
$$
Z(M)_{{e}{f}, {e}^\prime {f}^\prime}=
 \frac{\omega_{{e}}
\omega_{{e}^\prime}}{\omega^{m+1-\chi(\partial_+ M)/2}}
\sum_{{c}} \omega^2_{{c}}\;
Z(L_{{c}}\cup G^+_{{e}}\cup
G^-_{{e}^\prime}) \times $$
\be \label{z} \times
\frac{\omega_{{f}}\omega_{{f}^\prime}}
{\omega^{m+1-\chi(\partial_- M)/2}}
\sum_{{b}} \omega^2_{{b}}\;
Z(\bar{L}_{{b}}\cup \bar{G}^+_{{f}}\cup
\bar{G}^-_{{f}^\prime} )\, \ee
or 
\be \label{q}
Z(M)_{{e}{f}, \,{e}^\prime
{f}^\prime}=
\tau(M)_{{e}{e}^\prime}\;
\tau(-M)_{{f}^\prime f}\, .\ee

\vspace*{0.2cm}

{\bf Example:} Consider a solid torus $D^2\times S^1$.
Due to (\ref{z}) the corresponding state sum can be written as follows:
\vspace*{0.15cm}
\be\label{sa}
 Z_{ij}(D^2\times S^1)=\frac{\omega_i}{\omega^2} \sum_a \omega^2_a
Z(\hspace*{1.9cm}
\includegraphics{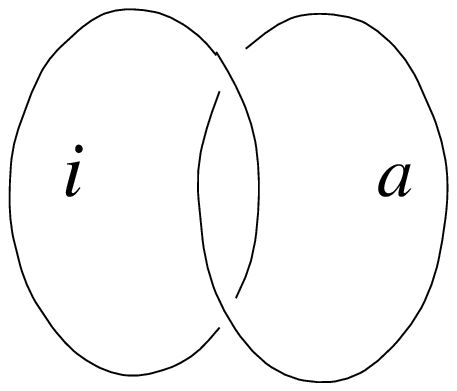} ) 
\;\frac{\omega_j}{\omega^2} \sum_b \omega^2_b
Z(\hspace*{1.9cm}
\includegraphics{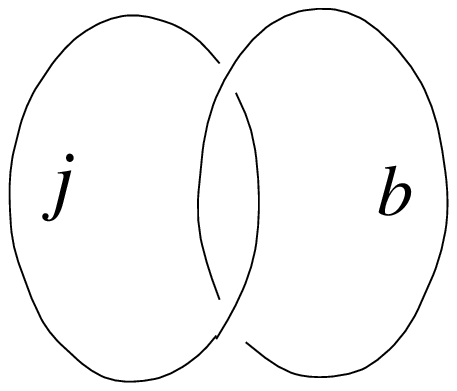} ). 
\ee
Recall that Euler characteristic of an empty set is equal to zero. 

We split the sums in 
(\ref{sa})  into the sums over even and odd colours, i.e.
$$
Z_{ij}(D^2\times S^1)= Z_{ij}(D^2\times S^1, s_0, 0)+
Z_{ij}(D^2\times S^1, s_1, 0)\, +$$
\be\label{sum} +
Z_{ij}(D^2\times S^1, s_0, h)+
Z_{ij}(D^2\times S^1, s_1, h), \ee
where
the first (resp. second) term 
corresponds in (\ref{sa})
to the case, when both $a$ and $b$ are even (resp. odd),
in the third term
$a$ is even and $b$ odd, and inversely in the forth term. Using
(\ref{circ1}) and (\ref{circ2}) one can calculate
$$ 
Z_{ij}(D^2\times S^1, s_0, 0)=\frac{1}{4}
(\delta_{i,0}+\delta_{i,\, {r-2}})(\delta_{j,0}+\delta_{j, \, r-2}),$$
$$ 
Z_{ij}(D^2\times S^1, s_1, 0)=\frac{1}{4}
(\delta_{i,0}-\delta_{i,\, {r-2}})(\delta_{j,0}-\delta_{j, \, r-2}),$$
$$ Z_{ij}(D^2\times S^1, s_0, h)=\frac{1}{4}
(\delta_{i,0}+\delta_{i,\, {r-2}})(\delta_{j,0}-\delta_{j, \, r-2}),$$
$$
Z_{ij}(D^2\times S^1, s_1, h)=\frac{1}{4}
(\delta_{i,0}-\delta_{i,\, {r-2}})(\delta_{j,0}+\delta_{j, \, r-2}).$$
Finally, we have
$$
Z_{ij}(D^2\times S^1)= \delta_{i,0}\;\delta_{j,0}\, .$$

\subsection[Refined Turaev--Viro TQFT]{Refined Turaev--Viro model}
In this section we 
refine the construction of [BD2].
 
\vspace*{0.2cm}
\centerline{\bf Definition of  invariants}
\vspace*{0.2cm}

We start by modifying   the notion of a
 standard handlebody.
  As before, consider the cylinder
 $\Sigma_g\times [0,1]$ with the graph ${\cal G}^g\in \Sigma_g$ and the graph
(\ref{basis}) on $-\Sigma_g$. 
Denote by $b_i$ a closed 1-dimensional subcomplex of the graph ${\cal
G}^g$, representing the $i^{\rm th}$ meridian of $\Sigma_g$, $1\leq
i\leq g$. We recall that the graph dual to ${\cal G}^g$ provides a
triangulation of $\Sigma_g$. 
Associate a ${\Bbb Z}_2$-number to the meridian $b_i$ of $\Sigma_g$,  
 $1\leq i\leq
g$, 
and denote by $\ma s$ a sequence of these numbers. Let $\ma h$
be a fixed subset 
of $\{b_i\}$.
These data define a standard handlebody $(M^+_g,{\ma s},{\ma h})$.

We say that  $({e}, {f})$ 
 is a special colouring 
of the graph $G^g\cup \bar{G}^g$, if the following conditions are
satisfied:
\begin{itemize}
\item
 colours $e_i$ and $f_i$, $1\leq i \leq g$, are even, if $b_i \notin  {\ma
h}$ and ${\ma s}_i=0$;  
\item
  colours $e_i$ and $f_i$, $1\leq i \leq g$, are odd, if $b_i \notin  {\ma
h}$ and ${\ma s}_i=1$;  
\item
 a colour $e_i$ is even and $f_i$ is odd,
 $1\leq i \leq g$,  if $b_i \in  {\ma
h}$ and ${\ma s}_i=0$;  
\item
 a colour $e_i$ is odd and $f_i$ is even,
 $1\leq i \leq g$,  if $b_i \in  {\ma
h}$ and ${\ma s}_i=1$.  
\end{itemize}
We denote the set of all special colourings by $E({\ma s}, {\ma h})$.
The  state sum for the standard handlebody 
 is given by the formula:
\be K_{{e}{f}} ({\ma s}, {\ma h}) =
\omega^{g-1}\omega_{{e}}
\omega_{{f}}
\sum_{{x}}
\prod^{3g-3}_{i=1} \frac
{\omega^2_{x_i}}{\omega^2}\; Z(M^+_g,
G^g_{{e}}\cup \bar{G}^g_f \cup m_{{x}} \cup
{\cal G}^g),\ee
and 
\be  K ({\ma s}, {\ma h}) =\oplus_{{e}\, {f}}
\; K_{{e} {f}} ({\ma s}, {\ma h}), \;\;\;\; (e,f)\in E({\ma s}, {\ma h}). \ee
This defines an inclusion
$$K ({\ma s}, {\ma h}): V_{\Sigma_g}({\ma s}, {\ma h})\to
V({\Sigma_g}) ,$$  where
$$   V_{\Sigma_g}({\ma s}, {\ma h}) =\oplus_{{e} \, {f}}
\;\{ V^L_g({{e}})\otimes V^R_g({f})\},\, 
 \;\;\;\; (e,f)\in E({\ma s}, {\ma h}). $$
 
The oppositely oriented  handlebody is given
 in the usual way as  the mirror
image of $(M^+_g, {\ma s}, {\ma h})$. The corresponding state sum 
$L({\ma s}, {\ma h})$  yields 
a projector
$$L ({\ma s}, {\ma h}):V({\Sigma_g})\to
 V_{\Sigma_g}({\ma s}, {\ma h}).$$
It is not difficult to verify by direct calculation  that
\be \label{id}
L({\ma s}^\prime ,{\ma h}^\prime) K({\ma s}, {\ma h})=\cases{0, & if
${\ma s}^\prime \neq {\ma s}$ and ${\ma h}^\prime \neq {\ma h}$;\cr
id_{V_{\Sigma_g}({\ma s}, {\ma h})}, & if
${\ma s}^\prime = {\ma s}$ and ${\ma h}^\prime = {\ma h}$.\cr}
\ee

By a parametrized triangulated surface $(\Sigma, {\ma s}, {\ma h})$ 
of genus $g$ we
mean a parametrized oriented
surface $\Sigma$ with triangulation $X$
provided with a sequence ${\ma
s}$ of ${\Bbb Z}_2$-numbers associated to the meridians $\phi(b_i)$,
$1\leq i\leq g$, and a fixed subset  $\phi({\ma
h})$
of the
 meridians, where $\phi: ({\cal G}^g)^\ast \to X$ is the parametrization 
 of $\Sigma$.
To the parametrized surface  $(\Sigma, {\ma s}, {\ma h})$ 
we associate  a vector space
$V_{\Sigma}({\ma s}, {\ma h})$, generated by  special colourings $E({\ma
s},{\ma h})$ of the graph $G^g\cup \bar{G}^g$.

\vspace*{0.2cm}

Consider  a 3-cobordism $(M,s,h)$ with parametrized boundary
$\partial M=
(-\partial_- M)\cup \partial_+  M$, where
$s$ is a spin structure on $M$ and $h\in H^1(M)$. Let us enumerate the
connected components of $\partial M$ by an index $j$, $1\leq j\leq
n$. Suppose that the first $l$ of them belong to $\partial_- M$ and
the remaining to $\partial_+ M$. Choose a sequence ${\ma s}_j$
of ${\Bbb Z}_2$-numbers and a set $ {\ma h}_j$
 on the $j^{\rm th}$ connected component $\Sigma_j$ of
$\partial M$, such that 
\be\label{param1}
({\ma
s}_j)_i=q_{s|_{\Sigma_j}}(\phi_j(b_i)),\;\; 1\leq i\leq g_j ,\ee
and ${\ma h}_j$ consists of the  meridians $b_i$,
such that $h$ is non-trivial on 
the homology class $[\phi_j(b_i)]\in H_1(M)$.
Here $\phi_j$ is the parametrization of $\Sigma_j$.

One glues  (along 
the  parametrizations)
  $(M^+_{g_j},{\ma s}_j, {\ma h}_j)$, $1\leq j\leq l$,
 and 
 $(M^-_{g_j},{\ma s}_j, {\ma h}_j)$, $l< j\leq n$, to the connected
components of $\partial_- M$ and $\partial_+ M$, respectively.
The resulting manifold can be represented by
surgery on $S^3$ with $n$ handlebodies removed and with a graph (given
by the image of  (\ref{basis}) under parametrization) sitting on the
boundary of each
handlebody  
(see [BD2] for more details). 
We set
$$
Z(M,s,h)_{{e}{f}, {
e}^\prime {f}^\prime
}=\omega^{{-\chi(\partial
M)}/{2}} \omega_{{e}} \omega_{{f}}
 \omega_{{e}^\prime} \omega_{{f}^\prime}
\sum_{{x} {y}{z}} \prod_{i\, j, k} \frac{\omega^2_{x_i}
\, \omega^2_{y_j}\, \omega^2_{z_k}}{\omega^2\, \omega^2\, \omega^2} \; $$
\be \label{Z1}
\sum_{{a} \,{b}\,a^\prime\, b^\prime}
Z(\tilde{S}^3, L_{{a}{b}}\cup m_{{z}}\cup
 G^+_{{e}{f}} \cup
G^-_{{e}^\prime{f}^\prime}\cup
m_{{x}}\cup m_{{y}})\, 
\prod^m_{i=1} S_{a_i \, a^\prime_i}\, S_{b_i \, b^\prime_i}\;
Z_{a^\prime_i\, b^\prime_i}(D^2\times S^1, s_i, h_i)\, ,
\ee
where 
$$  (e,f)\in E({\ma s}, {\ma h}), \;\;\;\;  (e^\prime,
f^\prime)\in E({\ma s}^\prime, {\ma h}^\prime),$$
$L$ is an $m$-component surgery link; $\tilde{S}^3$ is $S^3$  
with  neighborhoods of $L$, $G^+$ and $G^-$ removed; $L_{ab}\cup m_z$ is 
the coloured graph on the boundary of a  neighborhood of $L$
\footnotemark[4];
\footnotetext[4]{More precisely,   $L_{ab}\cup m_z=\cup^m_{i=1} (L_{a_i
b_i}\cup m_{z_i})$, where $L_{a_i b_i}$ consists of two ($a_i$- and
$b_i$-coloured) 
lines homotopic
to $L_i$,  where one of them  overcrosses  and the other one
undercrosses   meridian $m_{z_i}$.}
$ S_{ij}$ is an invariant of the  Hopf link (normalized by $\omega^{-1}$),
or equivalently,
 an element of MCG interchanging  cycles in the canonical
homology basis of a torus; $s_i$ and $h_i$ are the restrictions of $s$
and $h$  on the neighborhood of $L_i$. The state sums of a solid torus
with additional structures are listed in the  example of section 4.1,
where $s_0$ (resp. $s_1$) denotes the spin structure, which is (not)
bounding on $S^1$.

Taking into account that
$$\sum_{a^\prime} S_{a a^\prime}(\delta_{a^\prime ,\, 0}+
\delta_{a^\prime ,\, r-2})= \cases{ \omega^{-1}\omega^2_a, & if $ a$
is even\cr 0, & if $a$ is odd\cr}$$ 
$$\sum_{a^\prime} S_{a a^\prime}(\delta_{a^\prime ,\, 0}-
\delta_{a^\prime ,\,  r-2})= \cases{0, & if $a$
is even\cr \omega^{-1}{\omega^2_a}, & if $a$ is odd\cr}$$ 
and repeating the computation given in the proof of Theorem 2 in
[BD2], one obtains  that
\be \label{q2}
Z(M,s,h)_{{e}{f}, \,{
e}^\prime 
{f}^\prime}=
\tau(M,s)_{{e}{e}^\prime}\;
\tau(-M, s+h)_{{f}^\prime f}\, .\ee

As a result,  the operator $Z(M,s,h)$,
defined by  (\ref{Z1}), extends the Roberts' invariant to
an anomaly free non-degenerate TQFT.

\vspace*{0.2cm}
\centerline{\bf Gluing property}
\vspace*{0.2cm}



\begin{lem} If the  3-cobordism $(M,s,h)$ is obtained
from
$(M_1, s_1, h_1)$ and $(M_2, s_2, h_2)$ by gluing  along a
homeomorphism
 ${ f}: \partial_+ M_1\to \partial_- M_2$
 which preserves structure and commutes with parametrizations,
then  
$$\sum_{s,h}
Z(M,s,h)_{{e}{f}, \,{
e}^\prime 
{f}^\prime} = \sum_{e^{\prime \prime} f^{\prime\prime}}
Z(M_2, s_2, h_2)_{{e}{f}, \,{
e}^{\prime\prime} 
{f}^{\prime\prime}}\,
Z(M_1,s_1, h_1)_{e^{\prime\prime} f^{\prime \prime}, \,
{e}^\prime  {f}^\prime}\, ,$$
where the sum on the left hand side is taken over all $s$ and $h$, such
that $s|_{M_1}=s_1$, $s|_{M_2}=s_2$ and $h|_{M_1}=h_1$, $h|_{M_2}=h_2$.
\end{lem}

\vspace*{0.2cm}

\vspace*{0.2cm}
\centerline{ \bf Vector spaces associated to surfaces with structure}
\vspace*{0.2cm}

Due to (\ref{q2}),
for a closed  connected
surface $\Sigma$ with spin structure $\sigma$ and ${\rm h}\in H^1(\Sigma)$,
$$Z(\Sigma\times [0,1], \sigma\cup\sigma^\prime,{\rm h}\cup
{\rm h}^\prime)_{e f, {e}^\prime
{f}^\prime}
 = \cases{0, & if $\sigma\neq \sigma^\prime$ and
${\rm h}\neq {\rm h}^\prime$;  \cr
\tau^\sigma_{{e} {e}^\prime}
\tau^{\sigma+{\rm h}}_{{f}^\prime {f}}, & 
if $\sigma= \sigma^\prime$ and ${\rm h}={\rm h}^\prime$.  \cr
}
$$ 
Taking a direct sum over all special colourings we obtain an operator
$Z(\Sigma\times[0,1],\sigma,{\rm h})$. We define
the  vector space $V_{\Sigma}(\sigma,{\rm h})$ to be
 the support of this operator.
This vector space 
 is  associated by the spin TQFT of Turaev--Viro type to the
closed oriented connected surface $\Sigma$ provided with spin
structure $\sigma$ and first cohomology class h. 
Clearly,
$$ V({\Sigma})=\oplus_{\sigma,\, {\rm h}}
V_{\Sigma}(\sigma, {\rm h}),$$ 
$$ {\rm dim}V_\Sigma(\sigma, {\rm h})= 
{\rm dim}V_{(\Sigma,\, \sigma)}\,
{\rm dim}V_{(\Sigma,\, \sigma + {\rm h})}\, $$
and
$$Z(M,s,h): V_{\partial_- M}(s_-, {\rm h}_-)\to
V_{\partial_+ M}(s_+, {\rm h}_+)\, ,$$
where  
   $s_\pm=s|_{\partial_{\pm} M}$ and ${\rm h}_\pm=h|_{\partial_\pm M}$.

It follows from the results of Section 3.3, that
$$ Z( M)=\oplus_{s_\pm,\, {\rm h}_\pm }\sum_{s,h} Z(M,s,h),$$
where the sum is over $s$ and $h$, such that
 $s|_{\partial_{\pm} M}=s_\pm$ and $h|_{\partial_\pm M}={\rm h}_\pm$.
Moreover,
$$ Z( M,h )=\oplus_{s_\pm  }\sum_{s} Z(M,s,h)$$
is an invariant of a 3-cobordism $M$ with first cohomology class
$h$, which can be defined as follows (see [TV]):
Let us introduce a function  $a: I\to {\Bbb Z}_2$, such that
$$a(i)=i\pmod 2\, .$$
Then
for any admissible triple $(i,j,k)$
$$a(i)+a(j)+a(k)=0.$$
Therefore,
each colouring of a triangulated  3-manifold $M$ composed with $a$ is a
 1-cocycle of $M$. For any $h\in H^1(M)$, $Z(M,h)$ is equal to
 the Turaev-Viro invariant, where one sums  over all 
colourings  which induce cocycles representing $h$.  


\section{Concluding remarks}

In this article we restrict our attention to the case $r=0\pmod 4$,
because it corresponds to the invariants with the   richest topological
structure. The case $r=2\pmod 4$ can be treated by quite similar
methods, but it leads to invariants of 3-cobordisms  with a
first ${\Bbb Z}_2$-cohomology class only. For odd $r$ so far no 
refined invariants 
  are  known. 

It would be interesting to  find out
whether refined  quantum invariants  determined by additional topological
structures on 3-manifolds   could be defined for
higher quantum 
groups.
We leave this question for future investigation.

\vspace*{1cm}

{\Large \bf Acknowledgements}
\vspace*{1cm}

I would like to thank  Vladimir Turaev for many
stimulating discussions
and reading of the manuscript.

\end{document}